\documentclass[11pt,reqno]{article}

\usepackage[usenames]{color}
\usepackage{amssymb}
\usepackage{amsmath}
\usepackage{cases}
\usepackage{amsthm}
\usepackage{amsfonts}
\usepackage{amscd}
\usepackage{graphicx}

\usepackage[colorlinks=true,
linkcolor=webgreen,
filecolor=webbrown,
citecolor=webgreen]{hyperref}

\definecolor{webgreen}{rgb}{0,.5,0}
\definecolor{webbrown}{rgb}{.6,0,0}
\usepackage{color}
\usepackage{fullpage}
\usepackage{float}

\usepackage{cases}

\usepackage{graphics,amsmath,amssymb}
\usepackage{amsthm}
\usepackage{amsfonts}
\usepackage{latexsym}

\begin{document}

\theoremstyle{plain}
\newtheorem{theorem}{Theorem}
\newtheorem{remark}{Remark}
\newtheorem{lemma}{Lemma}

\newtheorem*{proposition}{Proposition}

\begin{center}
\vskip 1cm{\Large\bf A note on another approach on power sums} \\
%%\vskip .07in of powers of integers}
\vskip .2in \large Jos\'{e} Luis Cereceda \\
{\normalsize Collado Villalba, 28400 (Madrid), Spain} \\
\href{mailto:jl.cereceda@movistar.es}{\normalsize{\tt jl.cereceda@movistar.es}}
\end{center}

\begin{abstract}
In this note, we first review the novel approach to power sums put forward recently by Muschielok in arXiv:2207.01935v1, which can be summarized by the formula $S_m^{(a)}(n) = \linebreak \sum_{k} c_{mk} \psi_k^{(a)}(n)$, where the $c_{mk}$'s are the expansion coefficients and where the basis functions $\psi_m^{(a)}(n)$ fulfil the recursive property $\psi_m^{(a+1)}(n)= \sum_{i=1}^n \psi_m^{(a)}(i)$. Then, we point out a number of supplementary facts concerning the said approach not contemplated explicitly in Muschielok's paper. In particular, we show that, for any given $m$, the values of the $c_{mk}$'s can be obtained by inverting a matrix involving only binomial coefficients. This may be compared with the original approach of Muschielok, where the values of the $c_{mk}$'s can be obtained by inverting a lower triangular matrix involving the Stirling numbers of the first kind. Also, we make a conjecture about the functional form of the coefficients $c_{m\, m-k}$.
\end{abstract}

\section{Introduction}

For integers $m,a \geq 0$ and $n \geq 1$, the $a$-fold summation $S_m^{(a)}(n)$ (or hyper-sum) over the first $n$ positive integers to the $m$-th powers is defined recursively as
\begin{equation*}
S_m^{(a)}(n) = \left\{
              \begin{array}{ll}
                n^m  & \text{if}\,\, a =0, \\
                \sum_{i=1}^n S_m^{(a-1)}(i) &  \text{if}\,\, a \geq 1.
              \end{array}
            \right.
\end{equation*}
Several methods to obtain explicit formulas for $S_m^{(a)}(n)$ have been described in the literature (see, e.g., \cite{bounebirat,cere,cere2,chen,inaba,kargin,knuth,laissaoui}). Recently, in a very interesting paper, C. Muschielok \cite{musch} has developed another approach on the iterated power sums $S_m^{(a)}(n)$. Briefly, Muschielok's procedure is as follows.

First, define the polynomial sequence
\begin{equation*}
\psi_m(n) = n + (m-1)(n-1)B_{m-1,n-1},
\end{equation*}
where the binomial coefficient $B_{a,b}$ is given by
\begin{equation*}
  B_{a,b} = \binom{a+b}{a} = \binom{a+b}{b}.
\end{equation*}
Note that $\psi_m(n)$ can alternatively be expressed as
\begin{equation*}
  \psi_m(n) = B_{1,n-1} + m(m-1)B_{m,n-2},
\end{equation*}
where it is understood that $B_{m,n-2} =0$ when $n=1$. Similarly to $S_m^{(a)}(n)$, we can define recursively the polynomial series of $a$-th order $\psi_m^{(a)}(n)$ with respect to the sequence $\psi_m(n)$ as
\begin{equation*}
\psi_m^{(a)}(n) = \left\{
              \begin{array}{ll}
                \psi_m(n)  & \text{if}\,\, a =0, \\
                \sum_{i=1}^n \psi_m^{(a-1)}(i) &  \text{if}\,\, a \geq 1.
              \end{array}
            \right.
\end{equation*}
Furthermore, using the property $B_{a,b} = \sum_{\beta =1}^{b} B_{a-1,\beta}$, it can be shown that the polynomials $\psi_m^{(a)}(n)$ have the closed form \cite[Lemma 3]{musch}
\begin{equation}\label{basis}
  \psi_m^{(a)}(n) = B_{a+1,n-1} + \frac{m(m-1)}{m+a} (n-1) B_{m+a-1,n-1}.
\end{equation}
Having established the $\psi_m^{(a)}(n)$'s, the next step of Muschielok's procedure is to express the monomial $n^m$ as a linear combination of the polynomial basis $\psi_k(n)$'s, that is,
\begin{equation*}
  n^m = \sum_{k=0}^m c_{m k} \psi_k (n),
\end{equation*}
for some certain coefficients $c_{m k}$. Hence, from the preceding equation, it immediately follows that
\begin{equation*}
  S_m^{(a)}(n) = \sum_{k=0}^m c_{m k} \psi_k^{(a)}(n).
\end{equation*}
In order to find the coefficients $c_{m k}$, one expands $\psi_m(n)$ in terms of powers of $n$, that is,
\begin{equation*}
  \psi_m(n) = \sum_{i=0}^m a_{m i} n^i.
\end{equation*}
As this point, it should be noted that, as shown in \cite{musch}, $a_{m,0} = a_{m,1} =0$. Therefore, in the following (as is done in \cite{musch}), and without loss of generality, we restrict ourselves to the case where $m \geq 2$. Clearly, the bases $\{ n^k \}_{k=2}^{m}$ and $\{ \psi_k(n) \}_{k=2}^{m}$ form a pair of dual bases, and thus we have
\begin{equation*}
  \sum_{i=2}^m  a_{mi} c_{ij} = \delta_{mj}, \quad j = 2,3,\ldots,m,
\end{equation*}
or, equivalently,
\begin{equation}\label{rec}
c_{mk} =
\left\{
  \begin{array}{ll}
    1/a_{mm} = (m-2)!, & \hbox{if $k =m$;} \\
    -(m-2)! \sum_{l= k}^{m-1} a_{ml} c_{lk}, & \hbox{if $2 \leq k \leq m-1$;} \\
    0, & \hbox{else.}
  \end{array}
\right.
\end{equation}
The recurrence equation \eqref{rec} enables one to obtain the coefficients $c_{mk}$ in terms of the $a_{ml}$'s and the earlier coefficients $c_{lk}$. Alternatively, as pointed out in \cite{musch}, we can build the truncated matrix $A_m = (a_{ml})$ and calculate its inverse. As we shall presently see, the $a_{ml}$'s can be expressed in terms of the Stirling numbers of the first kind. For completeness, below we write down the matrix $C_{10} = (c_{mk})_{2 \leq m,k \leq 10}$ containing the coefficients $c_{mk}$ for $m =2,3,\ldots,10$, and where, for each involved $m$, the index $k$ ranges from $2$ to $m$:
\begin{equation*}
C_{10} =
\left(
  \begin{array}{ccccccccc}
    1 &   &   &   &   &   &   &   &   \\
    0 & 1 &   &   &   &   &   &   &    \\[-1mm]
    1 & -2 & 2 &   &   &   &  \text{\huge0} &   &    \\
    0 & 5 & -10 & 6 &   &   &   &   &    \\
    1 & -10 & 40 & -54 & 24 &   &   &   &  \\
    0 &  21 & -140 & 336 & -336 & 120 &   &   &   \\
    1 & -42 & 462 & -1764 & 3024 & -2400 & 720 &   &   \\
    0 & 85 & -1470 & 8442 & -22176 & 29520 & -19440 & 5040 &  \\
    1 & -170 & 4580 & -38178 & 144648 & -288000 & 313200 & -176400 & 40320  \\
  \end{array}
\right).
\end{equation*}
In this way, provided with the coefficients $c_{mk}$, we can obtain $S_m^{(a)}(n)$ from
\begin{equation}\label{proc}
  S_m^{(a)}(n) = \sum_{k=2}^m c_{mk} \psi_k^{(a)}(n),  \quad m \geq 2,
\end{equation}
with the basis functions $\psi_k^{(a)}(n)$ being given explicitly in \eqref{basis}.

\begin{remark}
The expansion coefficients $c_{mk}$ correspond to the sequence \href{https://oeis.org/A355570}{\emph{A355570}} in OEIS (The On-Line Encyclopedia of Integer Sequences).
\end{remark}

In the next section, we present a number of additional facts concerning Muschielok's approach not appearing in \cite{musch}. In particular, we show that, for any given $m$, the numerical values of the coefficients $c_{m2}, c_{m3}, \ldots, c_{mm}$ in the summation formula \eqref{proc} can be obtained by inverting a matrix involving only binomial coefficients (see equation \eqref{proc4} below). Moreover, we make a conjecture about the functional form of the coefficients $c_{m\, m-k}$, with $k \geq 0$ and $m \geq k+2$.

\enlargethispage{3mm}

\section{Facts and Conjecture}

In what follows, we point out several facts concerning Muschielok's summation procedure, along with a conjecture.

\vspace{3mm}
{\bf Fact 1.} When $a=0$, equation \eqref{proc} becomes $n^m = \sum_{k=2}^m c_{mk} \psi_k(n)$. Since $\psi_k(1) =1$, it immediately follows that
\begin{equation}\label{fact1}
  \sum_{k=2}^m c_{mk} =1, \quad\text{for all $m \geq 2$}.
\end{equation}
Property \eqref{fact1} can be readily checked for each of the rows in the matrix $C_{10}$ above.

\vspace{3mm}
{\bf Fact 2.} When $m=2$, from \eqref{proc} and \eqref{basis} we obtain
\begin{align*}
S_2^{(a)}(n) & = c_{22} \psi_2^{(a)}(n) \\[-1mm]
& = B_{a+1,n-1} + \frac{2}{a+2} (n-1) B_{a+1,n-1} \\[-1mm]
& = \frac{2n+a}{a+2} B_{a+1,n-1} = \frac{2n+a}{a+2} \binom{n+a}{a+1} = \frac{2n+a}{a+2} S_1^{(a)}(n),
\end{align*}
in accordance with the result for $S_2^{(a)}(n)$ given in \cite[p. 281]{knuth}.

\vspace{3mm}
{\bf Fact 3.} For $m\geq 2$ and $l =2,3,\ldots,m$, the coefficients $a_{ml}$ are given explicitly by
\begin{equation}\label{cofa1}
  a_{ml} = \frac{1}{(m-2)!} \left( \genfrac{[}{]}{0pt}{}{m-1}{l-1} - \genfrac{[}{]}{0pt}{}{m-1}{l} \right),
\end{equation}
or, equivalently,
\begin{equation}\label{cofa2}
  a_{ml} = \frac{1}{(m-2)!} \left( \genfrac{[}{]}{0pt}{}{m}{l} - m \genfrac{[}{]}{0pt}{}{m-1}{l} \right),
\end{equation}
where $\genfrac{[}{]}{0pt}{}{m}{l}$ denotes the (unsigned) Stirling numbers of the first kind. Using either \eqref{cofa1} or \eqref{cofa2}, we can construct the matrix $A_m = (a_{ml})$. For example, the matrix $A_{10} = (a_{ml})_{2 \leq m,l \leq 10}$ whose inverse is $C_{10}$, is given by
\begin{equation*}
A_{10} =
\left(
  \begin{array}{ccccccccc}
    1 &   &   &   &   &   &   &   &  \\
    0 & 1 &   &   &   &   &   &   &  \\
    -\frac{1}{2} & 1 & \frac{1}{2} &   &   &   &\text{\huge0}  &   &   \\[2mm]
    -\frac{5}{6} & \frac{5}{6} & \frac{5}{6} & \frac{1}{6} &   &   &   &   &  \\[2mm]
    -\frac{13}{12} & \frac{5}{8} & \frac{25}{24} & \frac{3}{8} & \frac{1}{24} &   &   &   &  \\[2mm]
    -\frac{77}{60} & \frac{49}{120} & \frac{7}{6} & \frac{7}{12} & \frac{7}{60} & \frac{1}{120} &   &   &  \\[2mm]
    -\frac{29}{20} & \frac{7}{36} & \frac{889}{720} & \frac{7}{9} & \frac{77}{360} & \frac{1}{36} & \frac{1}{720} &   &  \\[2mm]
    -\frac{223}{140} & -\frac{4}{315} & \frac{101}{80} & \frac{229}{240} & \frac{13}{40} & \frac{7}{120} & \frac{3}{560} & \frac{1}{5040} &  \\[2mm]
    -\frac{481}{280} & -\frac{61}{288} & \frac{1271}{1008} & \frac{427}{384} & \frac{853}{1920} & \frac{19}{192} & \frac{17}{1344} & \frac{1}{1152} & \frac{1}{40320} \\
  \end{array}
\right).
\end{equation*}
Let us further note, incidentally, that the elements of the first column of $A_m$ are given by $a_{m2} = 1- H_{m-2}$, $m\geq 2$, where $H_m$ denotes the $m$-th harmonic number.

\vspace{3mm}
{\bf Fact 4.} By substituting the coefficient $a_{ml}$ in \eqref{rec} by its expression in either \eqref{cofa1} or \eqref{cofa2}, one can derive explicit formulas for the coefficients $c_{mm}, c_{m \, m-1}, c_{m \, m-2}, c_{m \, m-3}$, etc., in succession. It should be noticed, however, that the complexity of the calculation of $c_{m \, m-k}$ grows rapidly with $k$. Next, we quote the exact formula of $c_{m \, m-k}$ for $k=0,\ldots,7$:
\begin{align*}
c_{mm} & = (m-2)!, \quad m \geq 2, \\[2mm]
c_{m \, m-1} & = (m-3)! \, (3-m) \frac{1}{2}m,  \quad m \geq 3, \\[2mm]
c_{m \, m-2} & = (m-4)! \binom{m}{2} \left[ \frac{1}{4} m^2 - \frac{23}{12}m + \frac{46}{12} \right],
\quad m \geq 4, \\[2mm]
c_{m \, m-3} & = (m-5)! \, (5-m) \binom{m}{3} \left[ \frac{1}{8} m^2 - \frac{9}{8}m + \frac{11}{4} \right],
\quad m \geq 5, \\[2mm]
c_{m \, m-4} & = (m-6)! \binom{m}{4} \left[ \frac{1}{16} m^4 - \frac{11}{8}m^3 + \frac{553}{48}m^2
- \frac{1747}{40}m + \frac{1901}{30} \right], \quad m \geq 6, \\[2mm]
c_{m \, m-5} & = (m-7)! \, (7-m) \binom{m}{5} \left[ \frac{1}{32} m^4 - \frac{37}{48}m^3 + \frac{697}{96}m^2
- \frac{1489}{48}m + \frac{611}{12} \right], \quad m \geq 7, \\[2mm]
c_{m \, m-6} & = (m-8)! \binom{m}{6} \left[ \frac{1}{64} m^6 - \frac{43}{64}m^5 + \frac{775}{64}m^4
- \frac{67513}{576}m^3 \right. \\
& \qquad\qquad\qquad\qquad\qquad\qquad\quad\quad\,\,  \left. + \, \frac{1930}{3}m^2 - \frac{1916141}{1008}m +
\frac{198721}{84} \right],  \quad m \geq 8, \\[2mm]
c_{m \, m-7} & = (m-9)! \, (9-m) \binom{m}{7} \left[ \frac{1}{128} m^6 - \frac{47}{128}m^5 + \frac{2777}{384}m^4
- \frac{88093}{1152}m^3 \right. \\
& \qquad\qquad\qquad\qquad\qquad\qquad\quad\quad\,\,\,\, \left. + \, \frac{14669}{32}m^2 - \frac{425993}{288}m +
\frac{16083}{8} \right],  \quad \, m \geq 9.
\end{align*}

Motivated by the patterns exhibited by the above expressions for $c_{m\, m-k}$, $k=0,\ldots,7$, we propose the following conjecture regarding the functional form of the coefficients $c_{m\, m-k}$:

\vspace{3mm}
{\bf Conjecture 1.}
\begin{itemize}
  \item For all integers $k \geq 0$ and $m \geq 2k+2$, we have
  \begin{equation*}
    c_{m\, m-2k} = (m-(2k+2))! \binom{m}{2k} \sum_{j=0}^{2k} \Gamma_{m,j} m^j,
  \end{equation*}
  where the non-zero (rational) coefficients in the set $\{ \Gamma_{m,0}, \Gamma_{m,1}, \ldots, \Gamma_{m,2k} \}$ have alternating signs, with the leading coefficient $\Gamma_{m,2k} = 1/2^{2k}$ being positive. Furthermore, the polynomial $\sum_{j=0}^{2k} \Gamma_{m,j} m^j$ is always a (rational) positive number, and the coefficient $c_{m\, m-2k}$ is a positive integer for all $m \geq 2k+2$. In particular, for $m=2k+2$, we have $c_{2k+2,2} =1$ for all $k \geq 0$.

  \item For all integers $k \geq 0$ and $m \geq 2k+3$, we have
  \begin{equation*}
    c_{m\, m-(2k+1)} = (m-(2k+3))! \, (2k+3-m) \binom{m}{2k+1} \sum_{j=0}^{2k} \Upsilon_{m,j} m^j,
  \end{equation*}
  where the non-zero (rational) coefficients in the set $\{ \Upsilon_{m,0}, \Upsilon_{m,1}, \ldots, \Upsilon_{m,2k} \}$ have alternating signs, with the leading coefficient $\Upsilon_{m,2k} = 1/2^{2k+1}$ being positive. Furthermore, the polynomial $\sum_{j=0}^{2k} \Upsilon_{m,j} m^j$ is always a (rational) positive number, and the coefficient $c_{m\, m-(2k+1)}$ is a negative integer for all $m \geq 2k+4$, whereas, for $m =2k+3$, we have that $c_{2k+3,2} =0$ for all $k \geq 0$.
\end{itemize}

\vspace{3mm}
{\bf Fact 5.} By using \eqref{basis} and \eqref{proc}, and taking into account the property \eqref{fact1} and the definition of $B_{a,b}$, we can express $S_m^{(a)}(n)$ in the form
\begin{equation}\label{proc2}
  S_m^{(a)}(n) = \binom{n+a}{a+1} \left[ 1+ (n-1)(a+1) \sum_{k=2}^m c_{mk} \binom{a+k}{k}^{-1}
\binom{n+a+k-2}{k-2} \right],  \quad m \geq 2.
\end{equation}
In particular, when $a=1$, the above expression yields the following formula for the ordinary power sums $S_m^{(1)}(n) = 1^m + 2^m + \cdots + n^m$:
\begin{equation*}
  S_m^{(1)}(n) = \frac{1}{2}n(n+1) \left[ 1+ 2(n-1) \sum_{k=2}^m \frac{c_{mk}}{k+1}
\binom{n+k-1}{k-2} \right],  \quad m \geq 2.
\end{equation*}
Note that formula \eqref{proc2} gives us $S_m^{(a)}(n)$ as $S_1^{(a)}(n)$ times a polynomial in $n$ of degree $m-1$. Furthermore, formula \eqref{proc2} tells us that $S_m^{(a)}(1) =1$ for all integers $a \geq 0$ and $m \geq 2$. (By the way, $S_m^{(a)}(1) =1$ for all integers $a \geq 0$ and $m \geq 0$.)

\vspace{3mm}
{\bf Fact 6.} For $a=0$, the formula \eqref{proc2} reduces to
\begin{equation*}
  n^m = n \left[ 1 +(n-1) \sum_{k=2}^m c_{mk} \binom{n+k-2}{k-2} \right],
\end{equation*}
which can be written in the form
\begin{equation}\label{proc3}
  \sum_{k=2}^m c_{mk} \binom{n+k-2}{n} = \frac{n^{m-1} -1}{n-1},
\end{equation}
provided that $m\geq 2$ and $n \geq 2$. Letting successively $n=2,3,\ldots,m$ in \eqref{proc3} yields the following linear system of $m-1$ equations in the unknowns $c_{m2}, c_{m3}, \ldots, c_{mm}$:
\begin{align*}
& \sum_{k=2}^m \binom{k}{2} c_{mk} = 2^{m-1} -1, \\
& \sum_{k=2}^m \binom{k+1}{3} c_{mk} = \frac{3^{m-1} -1}{2}, \\
& \quad \vdots  \\
& \sum_{k=2}^m \binom{m+k-2}{m} c_{mk} = \frac{m^{m-1} -1}{m-1}.
\end{align*}
Hence, solving for the coefficients $c_{mk}$, we are left with the matrix equation
\begin{equation}\label{proc4}
  \left(
    \begin{array}{c}
      c_{m2} \\[2mm]
      c_{m3} \\[3mm]
      \vdots \\[1mm]
      c_{mm} \\
    \end{array}
  \right) =
  \left(
    \begin{array}{cccccc}
      \binom{2}{2} & \binom{3}{2} & \binom{4}{2} & \ldots & \binom{m-1}{2} & \binom{m}{2} \\[2mm]
      \binom{3}{3} & \binom{4}{3} & \binom{5}{3} & \ldots & \binom{m}{3} & \binom{m+1}{3} \\[2mm]
      \vdots & \vdots & \vdots & \cdots & \vdots & \vdots \\[2mm]
      \binom{m}{m} & \binom{m+1}{m} & \binom{m+2}{m} & \ldots & \binom{2m-3}{m} & \binom{2m-2}{m} \\
    \end{array}
  \right)^{-1}
  \left(
    \begin{array}{c}
      2^{m-1} -1 \\[2mm]
      \frac{3^{m-1} -1}{2} \\[3mm]
      \vdots \\[1mm]
      \frac{m^{m-1} -1}{m-1}\\
    \end{array}
  \right),
\end{equation}
from which we can obtain, for any given $m$, the corresponding values of the $c_{mk}$'s. For example, for $m=10$, the matrix equation \eqref{proc4} reads
\begin{multline*}
\left(
  \begin{array}{l}
    c_{10,2} \\
    c_{10,3} \\
    c_{10,4} \\
    c_{10,5} \\
    c_{10,6} \\
    c_{10,7} \\
    c_{10,8} \\
    c_{10,9} \\
    c_{10,10} \\
  \end{array}
\right) =
\left(
  \begin{array}{ccccccccc}
    1 & 3 & 6 & 10 & 15 & 21 & 28 & 36 & 45 \\
    1 & 4 & 10 & 20 & 35 & 56 & 84 & 120 & 165 \\
    1 & 5 & 15 & 35 & 70 & 126 & 210 & 330 & 495 \\
    1 & 6 & 21 & 56 & 126 & 252 & 462 & 792 & 1287 \\
    1 & 7 & 28 & 84 & 210 & 462 & 924 & 1716 & 3003 \\
    1 & 8 & 36 & 120 & 330 & 792 & 1716 & 3432 & 6435 \\
    1 & 9 & 45 & 165 & 495 & 1287 & 3003 & 6435 & 12870 \\
    1 & 10 & 55 & 220 & 715 & 2002 & 5005 & 11440 & 24310 \\
    1 & 11 & 66 & 286 & 1001 & 3003 & 8008 & 19448 & 43758 \\
  \end{array}
\right)^{-1}
\left(
  \begin{array}{c}
    511 \\
    9841 \\
    87381 \\
    488281 \\
    2015539 \\
    6725601 \\
    19173961 \\
    48427561 \\
    111111111 \\
  \end{array}
\right) \\
=
\left(
  \begin{array}{c}
    1 \\
    -170 \\
    4580 \\
    -38178 \\
    144648 \\
    -288000 \\
    313200 \\
    -176400 \\
    40320 \\
  \end{array}
\right),
\end{multline*}
thus retrieving the values of $c_{10,2}, c_{10,3}, \ldots, c_{10,10}$ appearing in the last row of the matrix $C_{10}$.

\vspace{3mm}
{\bf Fact 7.} Let us recall that $S_m^{(a)}(n)$ admits the following representation in terms of the Stirling numbers of the second kind $\genfrac{\{}{\}}{0pt}{}{m}{k}$:
\begin{equation*}
  S_m^{(a)}(n) = \sum_{k=1}^{m} k! \genfrac{\{}{\}}{0pt}{}{m}{k} \binom{n+a}{a+k}, \quad a \geq 0, \, m \geq 1,
\end{equation*}
(see \cite{inaba}). Therefore, writing the last formula as
\begin{equation*}
  S_m^{(a)}(n) = \binom{n+a}{a+1} \left[ 1 + \binom{n+a}{a+1}^{-1} \sum_{k=2}^{m} k! \genfrac{\{}{\}}{0pt}{}{m}{k} \binom{n+a}{a+k} \right],
\end{equation*}
and comparing it to \eqref{proc2}, we obtain the identity
\begin{equation*}
\sum_{k=2}^m c_{mk} \binom{a+k}{k}^{-1} \binom{n+a+k-2}{k-2} =
\frac{1}{n(n-1)} \binom{n+a}{a}^{-1} \sum_{k=2}^{m} k! \genfrac{\{}{\}}{0pt}{}{m}{k} \binom{n+a}{a+k},
\end{equation*}
which holds for any integers $a\geq 0$, $m \geq 2$, and $n \geq 2$. In particular, when $a =0$, the above identity reduces to
\begin{equation*}
  \sum_{k=2}^m c_{mk} \binom{n+k-2}{k-2} = \frac{1}{n(n-1)} \sum_{k=2}^m  k! \genfrac{\{}{\}}{0pt}{}{m}{k} \binom{n}{k},
\end{equation*}
which, of course, is equivalent to \eqref{proc3}.

\vspace{3mm}
{\bf Fact 8.} Starting from $\psi_m(n) = \sum_{i=0}^m a_{m i} n^i$, one can readily obtain the recurrence relation
\begin{equation*}
 S_m^{(a)}(n) = (m-2)! \left[  \psi_m^{(a)}(n) - \sum_{i=2}^{m-1} a_{mi} S_i^{(a)}(n) \right],
\end{equation*}
which applies to any integers $a \geq 0$ and $m \geq 2$, and where the summation on the right-hand side is zero if $m=2$. For the case $m=2$, as shown in Fact 2, we have $S_2^{(a)}(n) = \psi_2^{(a)}(n) = \frac{2n+a}{a+2} \binom{n+a}{a+1}$. When $m >2$, the above recurrence relation gives us $S_m^{(a)}(n)$ in terms of $\psi_m^{(a)}(n)$, the coefficients $a_{mi}$, and the earlier sums $S_i^{(a)}(n)$, $i =2,3,\ldots, m-1$.

\vspace{2mm}

\section*{Note added}

In a sequel to \cite{musch}, Muschielok evaluated sums of the form
\begin{equation*}
T_m^{\alpha} = \sum_{k=2}^m c_{mk} k^{\alpha},
\end{equation*}
for integer $\alpha \geq 1$. In particular, he showed that $T_m^1 =  \sum_{k=2}^m c_{mk} k =m$ (\cite[Equation (35)]{musch2}). Moreover, he obtained the relation (\cite[Equation (48)]{musch2})
\begin{equation*}
n^m -n + m = \sum_{l=2}^n a_{nl} T_m^{l},
\end{equation*}
which we write as
\begin{equation}\label{seq1}
n^m -n = - T_m^1 + \sum_{l=2}^n a_{nl} T_m^{l},
\end{equation}
with $n \geq 2$. On the other hand, it is to be noted that \eqref{proc3} can be expressed in the form
\begin{equation*}
n^m -n = \frac{1}{(n-2)!} \sum_{k=2}^m c_{mk} (k-1)^{\overline{n}},
\end{equation*}
where $n \geq 2$ and $(k-1)^{\overline{n}} = (k-1)k\cdots (k+n-2)$. Hence, using that $(k-1)^{\overline{n}} = \sum_{j=0}^n \genfrac{[}{]}{0pt}{}{n}{j} (k-1)^j$, and applying the binomial theorem to $(k-1)^j$, it follows immediately that
\begin{equation}\label{seq2}
n^m - n = \sum_{l=0}^n A_{nl} T_m^l,
\end{equation}
where
\begin{equation*}
A_{nl} = \frac{1}{(n-2)!} \sum_{j=l}^n (-1)^{j-l} \binom{j}{l} \genfrac{[}{]}{0pt}{}{n}{j}.
\end{equation*}
Therefore, comparing \eqref{seq1} and \eqref{seq2}, we are led to conclude that $A_{n0} =0$, $A_{n1} = -1$, and, for $l \geq 2$, $A_{nl} = a_{nl}$. Renaming $n$ as $m$, the latter means that
\begin{equation}\label{seq3}
a_{ml} =  \frac{1}{(m-2)!} \sum_{j=l}^m (-1)^{j-l} \binom{j}{l} \genfrac{[}{]}{0pt}{}{m}{j}.
\end{equation}
Incidentally, in view of \eqref{cofa1} and \eqref{seq3}, we deduce the identity
\begin{equation*}
 \sum_{j=l}^m (-1)^{j-l} \binom{j}{l} \genfrac{[}{]}{0pt}{}{m}{j} = \genfrac{[}{]}{0pt}{}{m-1}{l-1} - \genfrac{[}{]}{0pt}{}{m-1}{l},
 \quad\text{with}\,\,\, 2\leq l \leq m.
\end{equation*}

\vspace{5mm}

\noindent
{\bf Acknowledgment}. I would like to thank Christoph Muschielok for useful interchange of e-mails concerning his papers \cite{musch} and \cite{musch2}.

\end{document}